\documentclass{amsart}
\usepackage[all]{xy}
\usepackage[latin1]{inputenc}
\usepackage{hyperref}
\usepackage{wasysym}

\newcommand{\cK}{{\mathcal K} }

\newcommand{\cT}{{\mathcal T} }

\newcommand{\cX}{{\mathcal X} }

\newcommand{\wt}{\widetilde}
\newcommand{\pt}{\partial}
\def\ol#1{{\overline{#1}}}
\newtheorem{theorem} {Theorem} [section]
\newtheorem{definition}[theorem] {Definition}
\newtheorem{lemma}[theorem]  {Lemma}

\newtheorem*{theorem*}{Theorem}
\newtheorem{proposition}[theorem] {Proposition}
\newtheorem*{proposition*}{Proposition}
\newtheorem{corollary}[theorem] {Corollary}
\newtheorem*{corollary*} {Corollary}

\def\ke{K{\"a}h\-ler-Ein\-stein }
\def\ks{Ko\-dai\-ra-Spen\-cer }
\def\ka{K{\"a}h\-ler}

\def\wp{Weil-Pe\-ters\-son }

\def\tei{Teich\-mül\-ler }
\def\ii{\sqrt{-1}}
\def\ddb{\sqrt{-1}\partial\overline{\partial}}

\def\cinf{C^\infty}

\def\ab{{\alpha\ol\beta}}

\def\gba{{g^{\ol\beta\alpha}}}
\def\gab{{g_{\alpha\ol\beta}}}
\def\na{\nabla_}

\begin{document}

\title[Variation of geodesic length functions]{Variation of geodesic length functions\\ over Teichmüller space}
\author[R.\ Axelsson]{Reynir Axelsson}
\address{Háskóli Íslands, Dunhaga 5, IS-107 Reykjavík, Ísland}
\email{reynir@raunvis.hi.is}
\author[G. Schumacher]{Georg Schumacher}
\address{Fachbereich Mathematik und Informatik,
Philipps-Universit\"at Marburg, Lahnberge, Hans-Meerwein-Strasse, D-35032
Marburg,Germany}
\email{schumac@mathematik.uni-marburg.de}
\date{}

\begin{abstract}
In a family of compact, canonically polarized, complex manifolds equipped
with \ke metrics the first variation of the lengths of closed geodesics
was previously shown in \cite{a-s} to be the geodesic integral of the
harmonic \ks form. We compute the second variation. For one dimensional
fibers we arrive at a formula that only depends upon the harmonic Beltrami
differentials. As an application a new proof for the plurisubharmonicity
of the geodesic length function and its logarithm (with estimate) follows,
which also applies to the previously not known cases of \tei spaces of
weighted punctured Riemann surfaces, where the methods of Kleinian groups
are not available.
\end{abstract}

\maketitle

\section{Introduction}\label{sec:intro}
In the study of \tei spaces geodesic length functions play an important
role, in particular under the aspect of the theory of several complex
variables.

In \cite{ke} Kerckhoff showed that for a finite number of closed
geodesics, which fill up a Riemann surface, the sum of the geodesic length
functions provides a proper exhaustion of the corresponding \tei space. In
\cite{wo-nielsen} Wolpert proved that this function is actually convex
along \wp geodesics and plurisubharmonic.  Later it turned out that the
logarithm of a sum of geodesic length functions is plurisubharmonic as
well \cite{wo:reprise,wo-jdg}. In \cite{ye} Yeung constructed  a bounded
plurisubharmonic exhaustion function together with estimates. The Levi
form of the geodesic length functions also played an important role in
McMullen's proof of the Kähler hyperbolicity of the moduli space
\cite{cmm}.

We want to base our study of geodesic length functions solely upon the
hyperbolic geometry of Riemann surfaces and use the methods of Kähler
geometry. From this point of view it is desirable to express results in
terms of harmonic Beltrami differentials, which are to be considered as
harmonic \ks forms.

This approach avoids entirely methods involving Fuchsian groups. In
particular our results extend to cases where uniformization theory is not
available, such as \tei and moduli spaces of weighted punctured Riemann
surfaces, equipped with conical hyperbolic metrics.

Our methods originate from the study of \ke manifolds of negative
curvature, and the computations are done in this framework. Our result for
the first variation of the geodesic length function is surprisingly
simple. In dimension one it reads:

\begin{theorem*}[{\cite[Theorem 3.2]{a-s}}]
Given a holomorphic family of hyperbolic Riemann surfaces $\cX \to S$, the
first variation of the length in a family of closed geodesics $\gamma_s$
is a geodesic integral of the harmonic Beltrami differential
$$
A_i=A^z_{i\ol z} \frac{\pt}{\pt z} \ol{dz}
$$
associated to a complex tangent vector $\pt/\pt s^i$, namely
\begin{equation*}
\frac{\partial \ell(\gamma_s)}{\partial s^i}=
\frac{1}{2} \int_{\gamma_s} A_i.
\end{equation*}
\end{theorem*}
When dealing with tensors of higher order like curvature, which involve
second order derivatives of metric tensors, certain integral operators
arise in a natural way. In the context of automorphic forms the operator
$$
(\Box + 1)^{-1},
$$
where $\Box$ denotes the (complex) Laplacian, was extensively studied
(cf.\ \cite{elstr}), and Wolpert used it in \cite{wo-curv}. Later it
played a major role in Siu's study of \ke manifolds \cite{siu:canlift} and
also in \cite{sch:curv,sch-preprint}.

Its counterpart for geodesic integration rather than integration over the
whole Riemann surface is the operator
$$
(-\frac{D^2}{dt^2} + c)^{-1}, \quad c=1,2,
$$
where $D/dt$ denotes covariant differentiation along a geodesic. Our main
theorem is the following:
\begin{theorem*}
Let $f:\cX \to S$ be a holomorphic family of hyperbolic Riemann surfaces
together with a differentiable family of closed geodesics $\gamma_s$. Then
\begin{gather*}
\frac{\pt^2 \ell(\gamma_s)}{\pt s^i \pt s^\ol\jmath} = \frac{1}{2}\Big(
\int_{\gamma_s}\Big( (\Box + 1)^{-1}(A_iA_\ol\jmath) + \big(-
\frac{D^2}{dt^2} + 2)^{-1}(A_i) \cdot A_\ol\jmath    \Big) \\ \hspace{6cm}
+ \frac{1}{2\ell(\gamma_s)}\int_{\gamma_s} A_i \cdot\int_{\gamma_s}
A_\ol\jmath \Big) .
\end{gather*}
\end{theorem*}
Let $ G^{W\!P}_{i\ol\jmath}=\langle A_i,A_j\rangle_{WP}$ denote the
coefficient of the \wp metric, and denote by $P_1$ a certain positive
function depending on the diameter (precisely, a lower bound for the
resolvent kernel).
\begin{corollary*}
The following estimate holds for the second variation:
\begin{eqnarray*}
\frac{\pt^2 \ell (\gamma_s)}{\pt s^i\pt s^\ol\jmath} &\geq& \frac{\ell(\gamma_s)}{2}
P_1(d(\cX_s)) \cdot \langle A_i,A_j\rangle_{WP} + \frac{1}{2\ell(\gamma_s)}
\int_{\gamma_s} A_i \int_{\gamma_s} A_\ol\jmath,
\\
\frac{\pt^2\log \ell (\gamma_s)}{\pt s^i\pt s^\ol\jmath} &\geq& \frac{1}{2}
P_1(d(\cX_s)) \cdot \langle A_i,A_j\rangle_{WP} .
\end{eqnarray*}
In particular, $\log \ell(\gamma_s)$ is strictly plurisubharmonic.
\end{corollary*}
(Here ''$\geq$'' is used in the sense of Definition~\ref{de:herm}.)
Further applications are given in Section~\ref{se:2ndtei}. We have upper
estimates.P
\begin{corollary*}
Let $\dim S=1$. Denote by $\|A_s\|_0$ the maximum of the pointwise norm of
the harmonic Beltrami differential taken over the fiber $\cX_s$. Then
\begin{eqnarray*}
\frac{\pt^2 \ell(\gamma_s)}{\pt s \pt \ol s}& \leq & \ell(\gamma_s) \|A_s\|^2_0,\\
\frac{\pt^2 \log\ell(\gamma_s)}{\pt s \pt \ol s}& \leq & \frac{3}{4} \|A_s\|^2_0.
\end{eqnarray*}
\end{corollary*}

\section{Families of \ke Manifolds}
We compute the second variations of the geodesic length function in the
general setting of \ke manifolds of negative Ricci curvature.

A Kähler form on a complex manifold  $X$  of dimension $n$  will be
denoted by
$$
\omega_X = \ii \gab dz^\alpha \wedge dz^\ol\beta.
$$
We use the summation convention together with the $\nabla$-notation for
covariant derivatives. A $|$--symbol will denote an ordinary derivative.
Also, $\partial_\alpha$ and $\partial_\ol\beta$ will stand for
$\partial/\partial z^\alpha$ and $\partial/\partial z^\ol\beta$
respectively. The raising and lowering of indices is defined as usual. We
also use the semi-colon notation for covariant derivatives. For the Ricci
tensor $R_{\alpha\ol\beta}$ on $X$ we use the sign convention
\begin{equation}\label{eq:ke}
R_\ab = - \log(g(z))_{|\ab},
\end{equation}
where $g(z)=\det(\gab(z))$.

Let $\{\cX_s\}_{s\in S}$ be a holomorphic family of canonically polarized
compact complex manifolds parameterized by a (connected) complex space
$S$. It is given by a proper, smooth, holomorphic mapping $f:\cX \to S$
such that $\cX_s=f^{-1}(s)$ for all $s\in S$. For simplicity we will
assume that the base $S$ is smooth, although our results can also be given
a meaning for possibly non-reduced singular base spaces.

Local coordinates on $S$ will be denoted by $s^i$, $i=1,\ldots, N$. We use
these as local coordinates on the total space $\cX$ together with further
local coordinates $z^\alpha$, $\alpha=1,\ldots, n$, where $n$ is the fiber
dimension, satisfying $f(z,s)=s$.

The fibers $\cX_s$ are equipped with  \ke forms
$$
\omega_{\cX_s} = \ii g_\ab(z,s) dz^\alpha\wedge dz^\ol\beta
$$
of constant negative curvature $-1$  depending smoothly upon the parameter
$s$. We write $g(z,s)= \det (g_\ab(z,s))$ for the family $f: \cX \to S$,
where $R_\ab(z,s) = - g_\ab(z,s)$.

We consider the real $(1,1)$-form
\begin{equation}\label{eq:omX}
\omega_\cX= \ddb \log g(z,s)
\end{equation}
on the total space $\cX$. The fiberwise \ke equation \eqref{eq:ke} implies
that
$$
\omega_\cX|\cX_s = \omega_{\cX_s}
$$
for all $s\in S$. In particular $\omega_\cX$, restricted to any fiber, is
positive definite. The following fact is known:
\begin{theorem*}[\cite{sch-preprint}]
Let $f:\cX \to S$ be nowhere infinitesimally trivial. Then $\omega_\cX$ is
a \ka\ form on the total space.
\end{theorem*}
Let
$$
\rho :T_{s}S \to H^1(\cX_s, \cT_{\cX_s})
$$
be the \ks map for the deformation $f:\cX \to S$ at a point $s\in S$.

The \ke metric $\omega_{\cX_s}$ on $\cX_s$ induces a natural inner product on the space $H^1(\cX_s, \cT_{\cX_s})$ of
infinitesimal deformations of $\cX_s$ and thus on $T_{s}S$ via the map $\rho$; this is the  {\it \wp}Hermitian inner product on
$T_{s}S$.
Namely, given tangent vectors $u,v \in T_{s}S$, we denote by $A_u=
A_{u\ol\beta}^\alpha
\partial_\alpha dz^\ol\beta$ and $A_v$ the harmonic representatives of
$\rho(u)$ and $\rho(v)$ respectively. Then the inner product of $u$ and
$v$ equals
$$
\langle u, v \rangle_{W\!P} = \int_{\cX_s} A_{u\ol\beta}^\alpha A_{\ol v
\gamma}^\ol\delta g_{\alpha\ol\delta}g^{\ol\beta\gamma} g \, dV,
$$
where $A_\ol v$ denotes the adjoint (conjugate) tensor of $A_v$, and $g\,
dV$ the volume element.

We note that the \wp inner product is positive definite at a given point of the
base, if the induced deformation is effective.

We set $A_j= A_{\partial/\partial s^j}$. Then the \wp form on $S$ equals
$$
\omega^{W\!P} = \ii G_{i\ol\jmath}^{W\!P}(s) ds^i \wedge ds^\ol\jmath,
$$
where we use the notation
$$
G_{i\ol\jmath}^{W\!P}(s)= \big\langle\partial /\partial s^i, \partial
/\partial s^j\big\rangle_{W\!P} = \int_{\cX_s} A_{i\ol\beta}^\alpha A_{\ol
\jmath \gamma}^\ol\delta g_{\alpha\ol\delta}g^{\ol\beta\gamma} g \,
dV.
$$
The short exact sequence
$$
0 \to \cT_{\cX/S} \to \cT_{\cX} \to f^*\cT_S \to 0
$$
induces the \ks map via the edge homomorphism for direct images.  A lift
of a tangent vector $\partial/\partial s^i$ at a point $s$ of $S$ is a
differentiable vector field on $\cX_s$ with values in $\cT_\cX $. It has
the form
$$
\partial/ \partial s^i + b_i^\alpha \partial_\alpha.
$$
Its exterior $\ol\partial$-derivative $B_{i\ol\beta}^\alpha
\partial_\alpha dz^\ol\beta$, where $B_{i\ol\beta}^\alpha=\na\ol\beta b_{i}^\alpha
$, is interpreted as a $\ol\partial$-closed $(0,1)$-form on $\cX_s$ with
values in the tangent bundle of $\cX_s$. Its cohomology class
\begin{equation}\label{eq:ks}
\rho(\partial/\partial s^i)= [B_{i\ol\beta}^\alpha \partial_\alpha
dz^\ol\beta] \in H^1(\cX_s,\cT_{\cX_s}).
\end{equation}
equals the obstruction against the existence of a holomorphic lift of the
given tangent vector, i.e.\ the infinitesimal triviality of the
deformation in the direction of the tangent vector.

We now introduce notations that will be used in the rest of the paper.

The {\em horizontal lift} of $\pt/\pt {s^i}$, i.e.\ the lift that is
perpendicular to the fibers, will be denoted by
\begin{equation}\label{eq:lift}
 v_i = \partial/\partial s^i + a_i^\alpha \partial_\alpha.
\end{equation}
Note that the quantities $a_i^\alpha$ are in general not tensors. It
follows from the definition that
\begin{equation}\label{eq:liftcomp}
a_i^\alpha= -\gba g_{i\ol\beta}.
\end{equation}
We set
\begin{equation}
A_{i\ol\beta}^\alpha = \na\ol\beta a_{i}^\alpha.
\end{equation}
The following properties of the tensors $A_{i\ol\beta}^\alpha$ are known
(cf.\ \cite{sch:curv,sch:teich}) and will be used in the sequel:

\begin{proposition}\label{pr:harm}
The horizontal lifts of tangent vectors with respect to $\omega_\cX$
induce the harmonic representatives of \ks classes in the sense that
$A_{i\ol\beta}^\alpha \partial_\alpha dz^\ol\beta$ is the harmonic
representative of $\rho(\partial/\partial s^i)$. The coefficients satisfy
the following properties
\begin{eqnarray}
\na\ol\delta A^\alpha_{i\ol\beta}&=&\na\ol\beta A^\alpha_{i\ol\delta},  \label{eq:harmi}\\
\na\gamma A_{i\ol\beta}^\alpha g^{\ol\beta\gamma}&=&0, \label{eq:harmii}\\
A_{i\ol\beta\ol\delta}&=&A_{i\ol\delta\ol\beta}. \label{eq:harmiii}
\end{eqnarray}
\end{proposition}
The conditions \eqref{eq:harmi} and \eqref{eq:harmii} above correspond to
harmonicity, whereas condition \eqref{eq:harmiii} reflects the
relationship with the metric tensor.

We use the notation $c^\ol\beta = \ol{c^\beta}$ for (locally defined)
tensors.

Later we will need the following fact:
\begin{lemma}
The partial derivatives of the Christoffel symbols with respect to the
base parameter satisfy the identities
\begin{eqnarray}
\Gamma^\alpha_{\gamma\sigma|s^i}&=& - a^\alpha_{i;\gamma\sigma}, \label{eq:Gammas} \\
\Gamma^\alpha_{\gamma\sigma|s^{\ol \jmath}}&=& - g^{\ol \beta \alpha}
 a_{\ol\jmath \gamma; \ol\beta \sigma}.\label{eq:Gammasq}
\end{eqnarray}
\end{lemma}

\section{Families of closed geodesics}
Let $(f: \cX \to S, \omega_\cX)$ be a family of \ke manifolds with
constant negative Ricci curvature $-1$, where $\omega_\cX$ is given by
\eqref{eq:omX}.

We denote by $\gamma_s$ a differentiable family of closed geodesics in the
fibers $\cX_s$, and by $\ell(s)$ the length of $\gamma_s$. In order to
compute first and second variations, it is sufficient to assume that $S$
is a disk in the complex plane centered at $0$ with coordinate $s$ (it is
even sufficient to assume that the embedding dimension equals one). The
general formulas follow from this case by polarization.

In local coordinates $(z,s)$ the closed geodesic curves $\gamma_s$ are
solutions of the differential equation
\begin{equation}\label{eq:geodesic}
\ddot u^\alpha(t,s) + \Gamma^\alpha_{\gamma\sigma}(u(t,s))\dot
u^\gamma(t,s)\dot u^\sigma(t,s) =0.
\end{equation}
The solution is unique up to an affine change of the parameter. In
particular we may prescribe any positive constant value of its speed
$$
\|\dot u(t,s)\|^2= {\gab(u(t,s),s)\dot u^\alpha(t,s)\dot
u^\ol\beta(t,s)}.
$$
For $s=0$ we choose $\|\dot u\|=1$, for the remaining values of $s$ the
value of $\|\dot u\|$ will be determined by the fact that the parameter
$t$ assumes values in the interval $[0, \ell_0]$, where $\ell_0$ is the
length of $\gamma_0$. Hence the family of geodesics is given by a map
\begin{equation}\label{eq:paraml0}
u: S \times [0,\ell_0] \to \cX
\end{equation}
such that $u \circ f $ is the projection onto the first factor. Now
\begin{equation}\label{eq:u_s}
u_*(\pt_s) = \pt_s + u^\alpha_s \pt_\alpha + u^\ol\beta_s \pt_\ol\beta
\end{equation}
with partial derivatives
$$
u^\alpha_s := u^\alpha_{|s}  \text{ and } u^\ol\beta_s := (\ol{u^\beta})_{|s}.
$$
Note that the $(1,0)$- and $(0,1)$-components $\pt_s + u^\alpha_s
\pt_\alpha$ and $u^\ol\beta_s \pt_\ol\beta$ of $u_*\pt_s$ are tensors
along the geodesics with values in $T_\cX$ and $T_{\cX/S}\subset T_\cX$
respectively. In a similar way the tensor
\begin{equation}\label{eq:u_t}
\dot u= u_*(d/dt) = \dot u^\alpha \pt_\alpha + \dot u^\ol\beta \pt_{\ol\beta}
\end{equation}
along the family of geodesics has a type decomposition. The difference of
two lifts of tangent vectors from the base is a tangent vector along the
geodesics (with values in the relative tangent bundle). For $\dim S =1$ we
have the horizontal lift
$$
v_s = \pt_s + a^\alpha\pt_\alpha.
$$

The difference of $u_*(\pt_s)$ and the horizontal lift has the components
\begin{eqnarray}
u^\alpha_s - a^\alpha_s &= &u^\alpha_s(s,t) - a^\alpha_s(u(s,t),s),\\
u^\ol\beta_s & =&  u^\ol\beta_s(s,t).
\end{eqnarray}

For any tensor along the geodesic $\gamma_s$ on a fiber $\cX_s$ we denote
by $D/dt$ the {\it covariant derivative along $\gamma_s$}. In particular
\begin{equation}\label{eq:Dtudot}
\frac{D}{dt}\dot u=0.
\end{equation}
Let $w^\alpha(t)\pt_\alpha$ be any vector field along $\gamma_s$. Then
\begin{equation}\label{eq:Ddt}
  \frac{D}{dt} w^\alpha(t) = \dot w^a(t) +
  \Gamma^\alpha_{\gamma\sigma}(u(t))w^\gamma(t)\dot u^\sigma(t).
\end{equation}
If $w^\alpha(t)$ is of the form $\wt w^\alpha(u(t))$, then \eqref{eq:Ddt}
implies
\begin{equation}\label{Ddttens}
  \frac{D}{dt} w^\alpha(t) = \wt w^\alpha(u(t))_{;\kappa}\dot u^\kappa(t) +
  \wt w^\alpha(u(t)_{;\ol \lambda} \dot u^\ol\lambda(t).
\end{equation}
Corresponding equations hold for tensors of type $(0,1)$ and contravariant
tensors.

\begin{lemma}\label{le:covderform}
We have
\begin{eqnarray}
  \frac{D}{dt}(u^\alpha_s - a^\alpha_s ) &=& \dot u^\alpha_s +
  \Gamma^\alpha_{\gamma\sigma} u^\gamma_s\dot u^\sigma - a^\alpha_{s;\gamma} \dot u^\gamma -
  a^\alpha_{s;\ol\beta} \dot u^\ol\beta \label{eq:D11}, \\
    \frac{D^2}{dt^2}(u^\alpha_s - a^\alpha_s )&=& \Gamma^\alpha_{\sigma\mu}\Gamma^\mu_{\kappa\gamma}
  \dot u^\kappa(\dot u^\sigma u^\gamma_s -\dot u^\gamma u^\sigma_s) + \label{eq:D21}\\ \nonumber
  &&  R^\alpha_{\; \sigma\gamma\ol\lambda}
  \dot u^\sigma\left(\dot u^\ol\lambda(u^\gamma_s -a^\gamma_s) - u^\ol\lambda_s \dot u^\gamma\right) -
  \\ \nonumber && 2 A^\alpha_{s \ol \delta; \gamma} \dot u^\gamma \dot u^\ol \delta -
  A^\alpha_{s \ol\delta;\ol\tau} \dot u^ \ol \delta \dot u^\ol \tau,\\
  \frac{D}{dt}(u^\ol\beta_s)&=& \dot u^\ol\beta_s +
  \Gamma^\ol\beta_{\ol\delta\ol\tau} u^\ol\delta_s \dot u^\ol\tau\label{eq:D12},\\
  \frac{D^2}{dt^2}(u^\ol\beta_s)&=& g^{\ol\beta\alpha} A_{s \ol\delta\ol\tau; \alpha } \dot u^\ol\delta \dot u^{\ol\tau}
  - R^\ol\beta_{\;\ol\delta\ol\tau\sigma} (u^\sigma_s -a^\sigma_s) \dot u^\ol\tau \dot u^\ol\delta + \label{eq:D22}
  \\ \nonumber &&   R^\ol\beta_{\; \ol\delta \ol\tau\sigma} \dot u^\sigma \dot u^\ol\tau u^\ol\delta_s +
  \Gamma^\ol\beta_{\ol \nu\ol\tau}\Gamma^\ol\nu_{\ol\delta\ol\lambda} \dot u^\ol\lambda
  (u^\ol\delta_s\dot u^\ol\tau - u^\ol\tau_s \dot u^\ol\delta ).
\end{eqnarray}
\end{lemma}
\begin{proof}
The equations \eqref{eq:D11} and \eqref{eq:D12} follow immediately from
the definition. The remaining proofs are rather computational: To prove
\eqref{eq:D21} we apply $D/dt$ to \eqref{eq:D11} and differentiate
\eqref{eq:geodesic} with respect to $s$. In this way we can eliminate
$\ddot u^\alpha_s$. We use \eqref{eq:Gammas}, and finally we have
\eqref{eq:D21}. Observe that we need to consider both ordinary and
covariant derivatives of Christoffel symbols. We prove \eqref{eq:D22} in
the same way.
\end{proof}

In order to describe the variation of the length of closed geodesics in a
family, we use the notion of integrating a tensor along a geodesic.
Exemplarily we define:
\begin{definition}
Let $C=C_{\ol\beta\ol\delta}$ be a tensor on the Kähler manifold $X$, and
$\gamma$ be a geodesic of length $\ell$, parameterized by $u(t)=
(u^1(t),\ldots,u^n(t))$, such that $\|\dot u(t)\|=1$. Then
$$
\int_\gamma C = \int_\gamma C_{\ol\beta\ol\delta}dz^\ol\beta
dz^\ol\delta := \int_0^\ell C_{\ol\beta\ol\delta}(u(t)) \dot
u^\ol\beta \dot u^\ol\delta dt.
$$
\end{definition}
For covariant tensors of order one this notation coincides with the
integration of a differential form along the curve $\gamma$. For contravariant
tensors the geodesic integral is defined after lowering indices with
respect to the metric tensor.

\section{First variation of the geodesic length function}
Given a holomorphic family of \ke manifolds with one dimensional base
space like in the previous section together with a differentiable family
of closed geodesics $\gamma_s$ with parametrization \eqref{eq:paraml0},
the length of these is equal to
$$
\ell(s)= \int_0^{\ell_0} \|\dot u(t,s)\| dt
$$
so that
\begin{equation}\label{eq:dls}
\left.\frac{d\ell(s)}{ds}\right|_{s=0}= \frac{1}{2} \int_0^{\ell_0}
\frac{d}{ds}\|\dot u(t,s)\| ^2 dt.
\end{equation}
We will compute
$$
\frac{d}{ds}\|\dot u(t,s)\| ^2 =\frac{d}{ds}\big( g_{\alpha\ol\beta}\dot u^\alpha\dot u^\ol\beta\big).
$$
We denote by $\langle \hspace{1.8mm},\hspace{1.2mm} \rangle_\cX $ the
inner product with respect to $\omega_\cX$.
\begin{lemma}
We have
\begin{equation}
\frac{d}{ds} \left(g_{\alpha\ol\beta}\dot u^\alpha\dot u^\ol\beta\right) -
\frac{d}{dt} \big\langle u_*\pt_s, \dot u\big\rangle_{\omega_\cX} =
A_{s\ol\beta\ol\delta}\dot u^\ol\beta \dot u^\ol\delta.
\end{equation}
\end{lemma}
In the computational {\em proof} one uses \eqref{eq:Dtudot},
\eqref{eq:Gammas}, \eqref{eq:D11}, and \eqref{eq:D12}.

An immediate consequence of the above Lemma is the following Theorem.

\begin{theorem}[{\cite[Theorem 3.2]{a-s}}]\label{th:varlength}
Given a holomorphic family of \ke manifolds with negative Ricci curvature,
the first variation of the length in a family of closed geodesics
$\gamma_s$ is the geodesic integral of the harmonic \ks tensors:
\begin{equation}\label{eq:dell}
\left.\frac{\partial \ell(s)}{\partial s^i}\right|_{s=s_0}\hspace{-1mm}=
\frac{1}{2} \int_{\gamma(s_0)} A_{i\ol\beta\ol\delta}dz^\ol\beta
dz^\ol\delta.
\end{equation}
\end{theorem}

\section{Second variation of the geodesic length function}
An important function is given by the inner product of harmonic lifts of
tangent vectors. In terms of local holomorphic coordinates $s^i$ on $S$
(or coordinates of a smooth ambient space of minimal dimension at a given
point of the base) we have:
\begin{definition}
Let $v_i$ be the horizontal lift of $\pt/\pt s^i$. We put
\begin{equation}
\varphi_{i\ol\jmath} = \langle v_i,v_j\rangle_{\cX},
\end{equation}
where the inner product is taken pointwise.
\end{definition}
We list basic properties of the function $\varphi_{i\ol\jmath}$ on $\cX$:
\begin{eqnarray}
\varphi_{i\ol\jmath} & = & g_{i\ol\jmath} - a^\alpha_i a^\ol\beta_\ol\jmath g_{\alpha\ol\beta}, \label{eq:proper1}\\
(\Box + 1)\varphi_{i\ol\jmath} & = & \langle A_{i\ol\beta\ol\delta}, A_{j \ol \lambda\ol\tau} \rangle=
A^\alpha_{i\ol\beta} A^\ol\beta_{\ol \jmath\alpha}, \label{eq:proper2}\\
\int_{\cX_s} \varphi_{i\ol\jmath}&=& G^{W\!P}_{i\ol \jmath}, \label{eq:proper3}\\
\omega^{n+1}_\cX & = & \ii \varphi_{i\ol\jmath} ds^i\wedge ds^\ol \jmath \wedge \omega^n_\cX.
\label{eq:proper4}
\end{eqnarray}
The first of these equalities follows from the definition. For the second
equality cf.\ \cite[Proposition 3]{sch-preprint}. The equation
\eqref{eq:proper3} follows from \eqref{eq:proper2}. The last equation
\eqref{eq:proper4} is Lemma~6 from \cite{sch-preprint}.

We will apply the following fact:
\begin{theorem*}[\cite{sch-preprint}]
The relative canonical bundle $\cK_{\cX/S}$ equipped with the hermitian
metric induced by the relative \ke forms is positive, i.e.\ the matrix
$(\varphi_{i\ol\jmath})$ is positive definite.
\end{theorem*}
The lower estimates for $(\varphi_{i\ol\jmath})$ from \cite{sch-preprint}
will be applied below.

Again, it is sufficient to do computations for a base space $S$ of
dimension one with coordinate $s$. By abuse of notation, we use $s$ and
$\ol s$ as indices instead of $i$ and $\ol\jmath$, where $i,j$ can only
take the value $1$.
\begin{lemma}
We have
$$
A_{s\ol\beta\ol\delta|\ol s}= -\varphi_{s\ol s;\ol\beta\ol\delta}
-A_{s\ol\tau\ol\beta;\ol\delta} a^\ol\tau_{\ol s} - A_{s\ol \tau\ol \beta} A^\ol\tau_{s\ol\delta}
- A_{s\ol\tau\ol\delta} A^\ol\tau_{\ol s \ol\beta}.
$$
\end{lemma}
\begin{proof}
We compute
\begin{gather*}
A_{s\ol\beta\ol\delta|\ol s}= \left(a_{s\ol\beta|\ol\delta} +
a_{s\ol\tau}\Gamma^\ol\tau_{\ol\beta\ol\delta}\right)_\ol s =
a_{s\ol\beta|\ol s;\ol \delta} + a_{s\ol\tau}
\Gamma^\ol\tau_{\ol\beta\ol\delta|\ol s}.
\end{gather*}
Now the claim follows from \eqref{eq:Gammas} and \eqref{eq:proper1}.
\end{proof}
From here we immediately obtain the following identity.
\begin{lemma}
We have
\begin{gather*}
\frac{\pt}{\pt\ol s}\big(A_{s\ol\beta\ol\delta}\dot u^\ol\beta\dot
u^\ol\delta\big) = (-\varphi_{s\ol s; \ol \beta\ol\delta} -
A_{s\ol\tau\ol\beta;\ol\delta} a^\ol\tau_\ol s - 2
A_{s\ol\tau\ol\beta}a^\ol\tau_{\ol s;\ol\delta}) \dot u^\ol\beta \dot
u^\ol\delta +\\ \hspace{30mm} A_{s\ol\beta\ol\delta|\ol\tau}u^\ol\tau_\ol
s \dot u^\ol\beta\dot u^\ol\delta +
A_{s\ol\beta\ol\delta;\gamma}u^\gamma_\ol s \dot u^\ol\beta \dot
u^\ol\delta + 2 A_{s\ol\beta\ol\delta} \dot u^\ol\beta_\ol s \dot u^\ol
\delta.
\end{gather*}
\end{lemma}
We need to eliminate mixed derivatives in the parameters $t$ and $s$. We
define a  function $\chi$ along the geodesics by the formula
\begin{equation}
\chi = \big\langle A^\kappa_{s \ol\lambda}\dot u^\ol\lambda \pt_\kappa, u_*(\pt_s)\big\rangle_{\omega_\cX} =
A_{s\ol\beta\ol\delta}(u^\ol\beta_\ol s -a^\ol\beta_\ol s)\dot u^\ol\delta
\end{equation}
and obtain
$$
\frac{d}{dt}\chi = \frac{D}{dt} \left(A_{s\ol\beta\ol\delta}\right)
(u^\ol\beta_\ol s -a^\ol\beta_\ol s)\dot u^\ol\delta
+A_{s\ol\beta\ol\delta}\frac{D}{dt} \left(u^\ol\beta_\ol s -a^\ol\beta_\ol s\right)\dot u^\ol\delta.
$$
A straightforward calculation using the identities \eqref{Ddttens} and
\eqref{eq:D11} shows that
\begin{gather}
\frac{\pt}{\pt\ol s}\big(A_{s\ol\beta\ol\delta}\dot u^\ol\beta\dot
u^\ol\delta\big) - 2\dot\chi + \frac{d}{dt}(\varphi_{s\ol s;\ol \beta}\dot
u^\ol\beta) = (\varphi_{s\ol s;\alpha\ol\beta} + 2 A_{s\ol\beta\ol\delta}
A^\ol\delta_{\ol s \alpha})\dot u^\alpha\dot u^\ol\beta
\\ \nonumber \hspace{40mm} - (A_{s\ol\beta\ol\delta;\ol\tau}\dot u^\ol\tau
+ A_{s\ol\beta\ol\delta;\gamma} \dot u^\gamma)(u^\ol\beta_\ol s -
a^\ol\beta_\ol s) \dot u^\ol\delta \\ \nonumber \hspace{40mm} +
A_{s\ol\beta\ol\delta;\gamma}\dot u^\ol\delta \left( u^\gamma_\ol s\dot
u^\ol\beta - (u^\ol \beta_\ol s - a^\ol\beta_\ol s  ) \dot u^\gamma
\right).
\end{gather}
This concludes the first part of the computation. Altogether we obtained:
\begin{proposition}\label{pr:2nd}
Let $\cX \to S$ be a holomorphic family of \ke manifolds of constant
negative Ricci curvature together with a differentiable family of closed
geodesics. Then the second variation of the geodesic length function
equals
\begin{gather}\label{eq:2nd}
\frac{\pt^2 \ell(s)}{\pt s\pt \ol s} = \frac{1}{2} \int_{\gamma_s}\Big(
(\varphi_{s\ol s;\alpha\ol\beta} + 2 A_{s\ol\beta\ol\delta}
A^\ol\delta_{\ol s \alpha})\dot u^\alpha\dot u^\ol\beta \hspace{30mm} \\
\nonumber \hspace{10mm}
 - (A_{s\ol\beta\ol\delta;\ol\tau}\dot u^\ol\tau
+ A_{s\ol\beta\ol\delta;\gamma} \dot u^\gamma)(u^\ol\beta_\ol s -
a^\ol\beta_\ol s) \dot u^\ol\delta   + A_{s\ol\beta\ol\delta;\gamma}\dot
u^\ol\delta \left( u^\gamma_\ol s\dot u^\ol\beta - (u^\ol \beta_\ol s -
a^\ol\beta_\ol s  ) \dot u^\gamma \right)\Big).
\end{gather}
\end{proposition}

\section{Second variation of the geodesic length function\\ on \tei
spaces}\label{se:2ndtei} From now on we assume that fibers of $f:\cX \to
S$ are one dimensional. We set $z=z^1$ and also use $z$ and $\ol z$ as
indices. The preceding formulas and the notation remain valid, if the
fibers are equipped with the hyperbolic metric of constant Ricci curvature
$-1$, i.e.\ on a fiber $\cX_s$ with coordinate function $z$ we have
$$
ds^2 = g(z,s) \ii dz\wedge \ol{dz}
$$
where $g(z,s)$ satisfies the equation
$$
g(z,s) = \frac{\pt^2 \log g(z,s)}{\pt z \ol{\pt z}}.
$$
Free homotopy classes of simple closed curves are represented by closed
geodesics $\gamma_s$ with parameterization $u(s,t)$, which depend in a
differentiable way upon the parameter $s$.

According to our general index convention we have $g= g_{z\ol z}$. Observe
that the harmonic \ks form
$$
A^z_{s\ol z} \pt_z \ol{dz}
$$
is exactly a harmonic Beltrami differential. Likewise
$$
A_{\ol s z z}
$$
defines a holomorphic quadratic differential. The statement of
Proposition~\ref{pr:2nd} now reads as follows:
\begin{proposition}\label{pr:2ndtei}
We have
\begin{gather}\label{eq:d2l}
\frac{\pt^2 \ell(\gamma_s)}{\pt s\pt \ol s} = \frac{1}{2}
\int_{\gamma_s}\Big( (\varphi_{s\ol s} + g^{\ol z z} A_{s\bar z\bar z}
A^\ol z_{\ol s z})
 + A_{s\bar z\bar z} \dot u^\ol z \frac{D}{dt}(u^\ol z_\ol s -
a^\ol z_\ol s)   \Big).
\end{gather}
\end{proposition}
\begin{proof}
The term that involves the function $\varphi$ can be interpreted as a
complex Laplacian and \eqref{eq:proper3} is applicable. We use $g\dot u
\dot{\ol u} =1$. The harmonicity of the \ks tensor is equivalent to
$$
A_{s  \bar z \bar z; z}=0
$$
so that the latter terms in \eqref{eq:2nd} vanish.
\end{proof}

\begin{theorem}\label{th:main}
Let $f:\cX \to S$ be a holomorphic family of hyperbolic Riemann surfaces
together with a differentiable family of closed geodesics $\gamma_s$. Then
\begin{gather}\label{eq:thm}
\frac{\pt^2 \ell(\gamma_s)}{\pt s^i \pt s^\ol\jmath} = \frac{1}{2}
\int_{\gamma_s}\Big( (\Box + 1)^{-1}(A_i\cdot A_\ol\jmath) + \big(-
\frac{D^2}{dt^2} + 2)^{-1}(A_i) \cdot A_\ol\jmath    \Big) \\ \hspace{6cm}
\nonumber +   \frac{1}{4\ell(\gamma_s)}\int_{\gamma_s} A_i
\cdot\int_{\gamma_s} A_\ol\jmath  .
\end{gather}
\end{theorem}
We prove the theorem in Section~\ref{se:pf}.

Next, we estimate the integrand in \eqref{eq:thm} from below:

\begin{definition}\label{de:herm}
Given any two Hermitian symmetric matrices $M_{i\ol \jmath}$ and $N_{i\ol
\jmath}$, we write $M_{i\ol \jmath}\geq N_{i\ol \jmath}$, if the
difference is a positive semi-definite matrix.
\end{definition}
\begin{corollary}\label{co:l2}
We have the inequality
$$
\frac{\pt^2 \ell(\gamma_s)}{\pt s^i \pt s^\ol\jmath} \geq
\frac{1}{2}\left( \int_{\gamma_s} (\Box + 1)^{-1}(A_i\cdot A_\ol\jmath) +
\frac{1}{\ell(\gamma_s)} \int_\gamma A_i \int_\gamma A_\ol\jmath
\right)
$$
in the sense of Definition~\ref{de:herm}. In particular the geodesic
length function is strictly plurisubharmonic.
\end{corollary}
Again we apply our main theorem from \cite{a-s}, which states that
$$
\frac{\pt\ell(\gamma_s)}{\pt s_i} = \frac{1}{2}\int_{\gamma_s} A_i
$$
and obtain the following statement.
\begin{corollary}\label{co:l2log}
The logarithm of the geodesic length function is strictly
plurisubharmonic: The inequality
$$
\frac{\pt^2\log \ell (\gamma_s)}{\pt s^i\pt s^\ol\jmath} \geq
\frac{1}{2 \ell(\gamma_s)}\int_{\gamma_s} (\Box +1)^{-1}(A_i\cdot A_\ol\jmath)
$$
holds in the sense of Definition~\ref{de:herm}.
\end{corollary}

A lower estimate for the functions $\varphi_{i\ol\jmath}=(\Box
+1)^{-1}(A_i\cdot A_\ol\jmath)$ is known:
\begin{proposition*}[cf.\ \cite{sch-preprint}]
There exists a positive function $P_1(d(\cX_s))$, which depends on the
diameter of $\cX_s$, such that for any solution
$$
(\Box + 1)\varphi = \chi,
$$
with $\chi\geq 0$ the inequality
$$
\varphi(z) \geq P_1(d(\cX_s)) \int_{\cX_s}\, \chi g\, dV
$$
holds for all $z \in \cX_s$.
\end{proposition*}
The above proposition implies the following estimate, which can be used
together with  Corollary~\ref{co:l2} and Corollary~\ref{co:l2log} to
obtain further inequalities:
$$
\int_{\gamma_s} (\Box +1)^{-1}(A_i\cdot A_\ol\jmath) \geq
\ell(\gamma_s) \cdot P_1(d(\cX_s)) \cdot G^{W\!P}_{i\ol\jmath}.
$$

\begin{lemma}[cf.\ {\cite[Lemma 3]{sch:framas}}]\label{le:sumpluri}
Let $\ell_j$ be positive functions on a complex manifold. Then the
following estimate of closed hermitian $(1,1)$-forms holds:
$$
\ddb \log\big(\sum_j \ell_j \big) \geq
\frac{1}{\sum_k \ell_k} \sum_j \big(\ell_j \ddb\log\ell_j\big).
$$
\end{lemma}

The above Lemma implies that estimates for the single geodesic length
functions carry over to any sum of such functions. Kerckhoff showed in
\cite{ke} that for a finite number of closed geodesics
$\gamma_1,\ldots,\gamma_m$, which fill up the Riemann surface the sum of
the geodesic length functions provides a proper exhaustion of the \tei
space. Wolpert proved in \cite{wo-nielsen} that the function
$(\ell(\gamma_1)+ \ldots + \ell(\gamma_m))^{1/2}$ is actually convex along
the \wp geodesics and $\log(\ell(\gamma_1)+ \ldots + \ell(\gamma_m))$ is
strictly plurisubharmonic (cf.\ \cite{wo:reprise,wo-jdg}).

Yeung constructs in \cite{ye} a bounded non-positive strictly
plurisubharmonic exhaustion function. His estimates of the second
variation of the geodesic length function follow from ours.

\begin{corollary}
The logarithm of any sum of geodesic length functions is strictly
plurisubharmonic with estimates given by Lemma~\ref{le:sumpluri}.
\end{corollary}

We conclude the section with an upper estimate, which we state for $\dim
S=1$.

\begin{corollary}\label{co:upper}
Let $\|A_s\|_0$ be the maximum of the pointwise norm of the harmonic
Beltrami differential taken over the fiber $\cX_s$. Then
\begin{eqnarray}
\frac{\pt^2 \ell(\gamma_s)}{\pt s \pt \ol s}& \leq & \ell(\gamma_s) \|A_s\|^2_0,\\
\frac{\pt^2 \log\ell(\gamma_s)}{\pt s \pt \ol s}& \leq & \frac{3}{4} \|A_s\|^2_0.
\end{eqnarray}
\end{corollary}
\begin{proof}
The maximum principle applied to the equation \eqref{eq:proper2} yields
that
$$
\varphi_{s\ol s}(z)\leq \|A_s\|^2_0.
$$
Furthermore,
$$
\int_{\gamma_s} \big(2- D^2/dt^2 \big)^{-1}(A_s)\cdot A_\ol s \leq
\frac{1}{2} \int_{\gamma_s} A_s \cdot A_\ol s \leq \frac{1}{2}\ell(\gamma_s) \|A_s\|^2_0 \;   ,
$$
and finally
$$
\big|\int_{\gamma_s}A_s\big|^2 \leq \ell(\gamma_s) \int_{\gamma_s} A_s \cdot A_\ol s.
$$
These estimates imply both inequalities.
\end{proof}

\section{Differential operators along closed geodesics}
When studying covariant differentiation along geodesics $u(t)$ on a fixed
Riemann surface, we observe that the obvious identities
$$
\frac{D}{dt} \dot u = 0 \qquad \text{and} \qquad \frac{D}{dt} g_{z\ol z}=0
$$
can be used to reduce covariant differentiation of tensors along a closed
geodesic to the (covariant) differentiation of functions. In our case all
functions will be of class $\cinf$. Hilbert space theory and regularity
theorems are available and need not explicitly be mentioned.
\begin{lemma}
The operator
$$
L=-\frac{D^2}{dt^2} + 1
$$
is invertible with bounded inverse.
\end{lemma}
Let $\lambda_\nu\geq 0$, $\lambda_0=0$ be the eigenvalues of
$-\frac{D^2}{dt^2}$. For any function $\psi$ we denote by
$$
\psi= \sum_{\nu\geq 0} \psi_\nu
$$
the eigenvector decomposition. An inverse of the operator $L-L^{-1}$ is
defined on the orthogonal complement $C$ of the kernel of $D^2/dt^2$
(which is also the kernel of $D/dt$) with values in the same complement.
\begin{lemma}\label{le:opM}
Let
$$
M=(L - L^{-1})^{-1}\circ\big( -\frac{D^2}{dt^2}\big).
$$
Then
$$
M(\psi)= \sum_{\nu>0} \big( 1 - \frac{1}{2 + \lambda_\nu}\big) \psi_\nu.
$$
In particular, $M= 1- (2 -D^2/dt^2)^{-1}$ on the complement $C$.
\end{lemma}

 This yields the following estimate.
\begin{proposition}\label{pr:opM}
For the geodesic integral we have
$$
0\leq \int_\gamma M(\psi)\ol \psi \leq \frac{1}{2} \Big(\int_\gamma |\psi|^2
- \frac{1}{\ell(\gamma)} \Big| \int_\gamma \psi  \Big|^2  \Big).
$$
\end{proposition}
\begin{proof}
By Lemma~\ref{le:opM} we have
$$
0\leq \int_\gamma M(\psi)\ol \psi= \sum_{\nu>0} \big( 1 - \frac{1}{2 + \lambda_\nu}\big)
\int_\gamma|\psi_\nu|^2 \leq \frac{1}{2} \int_\gamma \big(|\psi|^2 -
|\psi_0|^2\big).
$$
\end{proof}

\section{Proof of the main theorem}\label{se:pf}
We now prove Theorem~\ref{th:main}. Only the case of $\dim S=1$ is needed.

For one dimensional fibers \eqref{eq:D21} and \eqref{eq:D22} read
\begin{eqnarray}
\frac{D^2}{dt^2} (u^z_s -a^z_s) &=& (u^z_s -a^z_s) -g_{z\ol z} \dot u^z \dot u^z u^\ol z_s
-A^z_{s \ol z;\ol z} \dot u^\ol z \dot u^\ol z, \label{eq:D211}  \\
\frac{D^2}{dt^2} (u^z_\ol s) & = & -
g_{z\ol z} (u^\ol z_\ol s - a^\ol z_\ol s) \dot u^z \dot u^z + u^z_\ol s. \label{eq:D221}
\end{eqnarray}

We define auxiliary functions along the geodesics. Let
\begin{eqnarray*}
w &=& (u^z_s - a^z_s) \dot u^\ol z g_{z\ol z},\\
v &=& u^\ol z_s \dot u^z g_{z\ol z},\\
A &=& A_{s\ol z\, \ol z}(\dot u^\ol z)^2.
\end{eqnarray*}
We apply \eqref{eq:D211} and \eqref{eq:D221} and use the notation of the
preceding paragraph. The aim is to express the function $w$ in terms of
the \ks form. We have
\begin{eqnarray*}
Lw &=& v + \frac{D}{dt} A,\\
Lv &=& w.
\end{eqnarray*}
Since $(D/dt)A$ is orthogonal to the kernel of $D^2/dt^2$, we have
$$
\frac{D}{dt} A = (L-L^{-1})w.
$$
Now
\begin{gather*}
\int_{\gamma_s} w \frac{D}{dt}\ol A = \int_{\gamma_s}\Big(
(L-L^{-1})^{-1} \frac{D}{dt} (A) \Big) \cdot\frac{D}{dt}(\ol A) =\\
\hspace{30mm}-\int_{\gamma_s} \Big( (L-L^{-1})^{-1} \frac{D^2}{dt^2} (A)
\Big) \cdot\ol A = \int_{\gamma_s} M(A)\cdot\ol A.
\end{gather*}
Now Theorem~\ref{th:main} follows from Proposition~\ref{pr:2ndtei} and
Lemma~\ref{le:opM}. \qed

\section{Weighted punctured Riemann surfaces and conical metrics}\label{sec:cone}
In our previous paper \cite{a-s}, we discussed the first variation of the
geodesic length function for \tei spaces of weighted punctured Riemann
surfaces equipped with hyperbolic conical metrics. Using the extended
techniques in \cite{sch:conic} one can see that our results on second
variations and plurisubharmonicity hold true in the conical case (for
weights $\geq 1/2)$.

\end{document}